\documentstyle[12pt]{article}
\def \Z{\hbox{$Z\hskip -5.2pt Z$}}
\def\DD{\mbox{$\cal D$}}
\def\Z{\hbox{$Z\hskip -5.2pt Z$}}

\def\dis{\displaystyle}
\def\sc{\scriptstyle}

\def\cl{\centerline}
\def\ol{\overline}
\def\ul{\underline}
\def\rar{\rightarrow}
\def\Lar{\Leftarrow}
\def\Rar{\Rightarrow}
\def\Lra{\Leftrightarrow}

\def\bs{\backslash}
\def\vs{\vspace}
\def\a{\alpha}
\def\b{\beta}
\def\d{\delta}
\def\e{\epsilon}
\def\g{\gamma}
\def\l{\lambda}
\def\p{\partial}
\def\G{\Gamma}
\begin{document}
\cl {{\bf SIMPLE LIE COLOR ALGEBRAS OF WEYL TYPE${\sc\,}$}\footnote
{ AMS Subject Classification - Primary: 17B20, 17B65, 17B67,
17B68, 17B75,
16A68, 16A40.\\
\indent\hskip .3cm This work was supported by NSF of China,
National Educational Department of China, Jiangsu Educational
Committee, and Hundred Talents Program of Chinese Academy of
Sciences.}} {\small\vs{4pt} \par\cl{(To appear in {\it Israel J.
Math.})}\vs{4pt}\par \cl{Yucai Su\footnote { These authors were
partially supported by Academy of Mathematics and System Sciences
during their visit to this academy.}, \,\,\,\, Kaiming
Zhao\,\,\,and\,\,\, Linsheng Zhu$^2$}

\vs{8pt}
\par
{\bf Abstract.} For an $(\e,\G)$-color-commutative associative
algebra $A$ with an identity element over a field $F$ of
characteristic not 2, and for a color-commutative subalgebra $D$
of color-derivations of $A$, denote by $A[D]$ the associative
subalgebra of ${\rm End}(A)$ generated by $A$ (regarding as
operators on $A$ via left multiplication) and $D$. It is easily
proved that, as an associative algebra, $A[D]$ is $\G$-graded
simple if and only if $A$ is $\G$-graded $D$-simple. Suppose $A$
is $\G$-graded $D$-simple. Then,
 (a) $A[D]$ is a free left $A$-module;
 (b) as a Lie color algebra, the subquotient
     $[A[D],A[D]]/Z(A[D])\cap[A[D],A[D]]$ is simple
     (except one minor case),
where $Z(A[D])$ is the color center of $A[D]$. The structure of this
subquotient is explicitly described.
}\vskip 10pt \centerline{\bf \S1 Introduction}
\par
Lie color algebras are generalizations of  Lie superalgebras. Let
us start with the definition. Let $F$ be a field, $\G$   an
additive group. A {\it skew-symmetric bicharacter} of $\G$ is a
map
 $\e:\G\times\G\rar F^*$
satisfying
$$
\e(\l,\mu)=\e(\mu,\l)^{-1},\,
\e(\l,\mu+\nu)=\e(\l,\mu)\e(\l,\nu),\,
\forall\,\l,\mu,\nu\in\G.
\eqno(1.1)$$
It is clear that $\e(\l,0)=1$ for any $\l\in\G$.
Let $L=\oplus_{\l\in\G}L_{\l}$ be a $\G$-graded
$F$-vector space. For a nonzero homogeneous element $a\in L$, denote by
$\bar{a}$ the unique group element in $\G$ such that
$a\in L_{\bar{a}}$. We shall call $\bar{a}$ the {\it color} of $a$.
The $F$-bilinear map $[\cdot,\cdot]: L\times L\rar L$ is called a
Lie color bracket on $L$ if the following conditions are satisfied:
$$
\matrix{
 [a,b]=-\e(\bar{a},\bar{b})[b,a],\ \ \mbox{(skew symmetry)}
\vs{4pt}\hfill\cr
 [a,[b,c]]=[[a,b],c]+\e(\bar{a},\bar{b})[b,[a,c]],\ \
 \mbox{(Jocoby identity)}\cr
}
\eqno\matrix{(1.2)\vs{4pt}\cr(1.3)\cr}
$$
for all homogeneous elements $a,b,c\in L$. The algebra structure
$(L,[\cdot,\cdot])$ is called an {\it $(\e,\G)$-Lie color algebra}
or simply a {\it Lie color algebra}. If $\G=\Z/2\Z$ and
$\e(i,j)=(-1)^{ij},\forall\,i,j\in \Z/2\Z$, then $(\e,\G)$-Lie
color algebras are simply Lie superalgebras. For Lie color
algebras, we refer the reader to [P2] and the references there.
\par
Let $A=\oplus_{\l\in\G}A_{\l}$ be a $\G$-graded associative
$F$-algebra with an identity element 1, i.e.,
$A_{\l}A_{\mu}\subset A_{\l+\mu}$ for all $\l,\mu\in\G$. So $1\in
A_0$. We say that $A$ is {\it graded simple} if $A$ does not have
nontrivial $\G$-graded ideals. Denote by $H(A)$ all homogeneous
elements of $A$. If we define the bilinear product $[\cdot,\cdot]$
on $A$ by
$$
[x,y]=xy-\e(\bar{x},\bar{y})yx,\ \forall\,x,y\in H(A),
\eqno(1.4)$$ then $(A,[\cdot,\cdot])$ becomes a $(\e,\G)$-Lie
color algebra. We shall simply write $\e({x},{y})$ for
$\e(\bar{x},\bar{y})$ for $x,y\in H(A)$.
\par
An {\it Lie color ideal} $U$ of $A$ is a $\G$-graded vector space
$U$ of $A$ such that $[A,U]\subset U$. Sometimes it is called an
$(\e,\G)$-Lie ideal, or simply a color ideal. The {\it
$\e$-center} $Z(A)$ of $A$ is defined as
$$
Z=Z(A)=\{x\in A\,|\,[x,A]=0\}.
\eqno(1.5)$$
It is easy to see that $Z (A)$ is $\G$-graded. We say that $A$ is
{\it color-commutative} (or $\e$-color commutative) if $[A,A]=0$.
\par
Passman [P1,P2] proved that, for a color-commutative associative
algebra $A$ with an identity element over a field $F$, and for a
color-commutative color-derivation subalgebra $D$ of $A$, the Lie
color algebra (including the Lie algebra case) $AD=A\otimes D$ is
simple if and only if $A$ is graded $D$-simple and $AD$ acts
faithfully on $A$ and, ${\rm char\sc\,}F\ne2$ or ${\rm
dim}_FD\ge2$ or $D(A)=A$. In [SZ1, SZ2, Z],  (associative and Lie)
algebras of Weyl type were constructed and studied. In this paper
we shall study the color version of these algebras.
\par
For a color-commutative associative algebra $A$ with an identity
element over a field $F$ of characteristic not 2, and for a
color-commutative subalgebra $D$ of color-derivations of $A$,
denote by $A[D]$ the associative subalgebra of ${\rm End}(A)$
generated by $A$ (regarding as operators on $A$ via
multiplication) and $D$. It is easily proved that, as an
associative algebra, $A[D]$ is graded simple if and only if $A$ is
graded $D$-simple (see Theorem 2.3). Suppose $A$ is graded
$D$-simple. Then, (a) $A[D]$ is a free left $A$-module; (b) as a
Lie color algebra, the subquotient
  $[A[D],A[D]]/Z(A[D])\cap[A[D],A[D]]$ is simple
  (except one minor case),
where $Z(A[D])$ is the $\e$-center of $A[D]$.
The structure of this
subquotient is explicitly described (see Theorem 3.9). In many cases
$A[D]=[A[D],A[D]]$.

\vskip 10pt
\centerline{\bf \S2 Graded simple associative algebras of Weyl type}
\par

Throughout this work, we assume that  $A$ is $\G$-graded and $(\e,
\G)$-color commutative and the field $F$ is of characteristic not
=2. An $F$-linear transformation $\p:A\rar A$ is called a {\it
homogeneous color-derivation} of degree $\bar{\p}\in\G$ if
$$
\matrix{
 \p(a)\in A_{\bar{\p}+\bar{a}},\ \forall\,a\in H(A)\ {\rm and}
\vs{4pt}\hfill\cr
 \p(ab)=\p(a)b+\e({\p},{a})a\p(b),\ \forall\,a,b\in H(A).
\hfill\cr} \eqno(2.1)$$ Taking $a=b=1$ in (2.1), we obtain
$\p(1)=0$ and so $\p(c)=0$ for all $c\in F$. Denote ${\rm
Der}^\e(A)=\oplus_{\l\in\G}{\rm Der}_\l^\e(A)$, where ${\rm
Der}_\l^\e(A)$ is the $F$-vector space spanned by all homogeneous
color derivations of degree $\l\in\G$. As in the Lie algebras
case, it is easy to verify that ${\rm Der}^\e(A)$ becomes an $(\e,
\G)$-Lie color algebra under the Lie color bracket
$$
[\p,\p']=\p\p'-\e({\p},{\p'})\p'\p,\ \forall\,\p,\p'\in H({\rm
Der}^\e(A)), \eqno(2.2)$$ where $\p\p'$ is the composition of the
operators $\p$ and $\p'$.
\par
Let $D=\oplus_{\l\in\G}D_\l$ be a nonzero $\G$-graded
color-commutative Lie color subalgebra of ${\rm Der}^\e(A)$, i.e.,
$$
\p\p'=\e({\p},{\p'})\p'\p,\ \forall\,\p,\p'\in H(D). \eqno(2.3)$$
We call $A$ is {\it graded $D$-simple} if $A$ has no nontrivial
graded $D$-stable ideals (see [P2]). Set
$\G_+=\{\l\in\G\,|\,\e(\l,\l)=1\}$. Then by (1.1), $\G_+$ is a
subgroup of $\G$ with index $\le2$. Set
$\G_-=\{\l\in\G\,|\,\e(\l,\l)=-1\}$. For any graded subspace $B$
of $A$, we define $B_+=\oplus_{\l\in\G_+}B_\l$, then $B_+$ is
$\G$-graded. Similarly we can define $B_-$. Since
$\G=\G_+\cup\G_-$, it follows that $B=B_+\oplus B_-$. By (2.3), we
have $ \p^2=0\mbox{ if } \p\in H(D_-). $

\par
Fix a homogeneous basis $\{\p_i\,|\,i\in I\}$ of $D$, where ${I}$
is some index set. Fix a total ordering $<$ on ${I}$. Define $ J $
to be the set of all
 $\a=(\a_i\,|\,i\in I)\in\Z_+^{\,\,{ I}}$  such that $\a_i=0$ for all but a finite
number of $i\in I $ and $ \a_i=0\mbox{ or }1\mbox{ if
}\bar\p_i\in\G_-. $ For $\a\in J $, we define the {\it support} of
$\a$ by ${\rm supp}(\a)=\{i\in I \,|\,\a_i\ne0\}$ and define the
{\it level} of $\a$ by $l(\a)=|\a|=\sum_{i\in I}\a_i$. Define a
total ordering on $ J $ by
$$
\a<\b \Lra |\a|<|\b|,\mbox{ or }|\a|=|\b|\mbox{ and }\exists\,i\in
 I \mbox{ such that } \a_i<\b_i\mbox{ and
}\a_j=\b_j,\forall\,j\in  I ,j< i. \eqno(2.4)$$ Denote by $F[D]$
the $\G$-graded associative color-commutative algebra generated by
$D$, and denote by $A[D]$ the associative subalgebra of ${\rm
End}(A)$ generated by $A$ (regarding as operators on $A$ via left
multiplication) and $D$. For $x\in A,\a\in J $, we define by
$x\p^\a$ the operator acting on $A$, i.e.,
$$
(x\p^\a)(y)=x(\p_1^{\a_1}(\p_2^{\a_2}...(\p_r^{\a_r}(y))...)),
\forall\,y\in A, \eqno(2.5)$$ where
$\p^\a=\p_{i_1}^{\a_{i_1}}\p_{i_2}^{\a_{i_2}}
\cdots\p_{i_r}^{\a_{i_r}}$ with $i_1< i_2< \cdots<i_r$. We shall
call $\a$ the {\it degree} of $\p^{\a}$. Define $\e^+({\a,\b})$
for $\a,\b\in \G$ by $\p^\a\p^\b=\e^+({\a,\b})\p^{\a+\b}$, then
$$
\e^+({\a,\b})=\prod_{i,j\in{\cal
I}:\,i>j}\e({{\p_i}},{\p_j})^{\a_i\b_j},\, \forall\,\a,\b\in{ J}.
\eqno(2.6)$$ We shall also simply write $\e^+({x},{y})$ for
$\e^+(\bar{x},\bar{y})$ for $x,y\in H(A[D])$. From
$\e^+({\a,\b})\p^{\a+\b}=\p^\a\p^\b=\e({\a,\b})\p^{\b}\p^{\a}
=\e({\a,\b})\e^+({\b,\a})\p^{\a+\b}$, we see that
$$\e(\a,\b)=\e^+({\a,\b})\e^+({\b,\a})^{-1},\,
\forall\, \a,\b\in J , \eqno(2.7)$$ If $\a+\b\notin J $, we set
$\p^{\a+\b}=0$. For any $\a\in J $, set
$$ J (\a)=\{\g\in J \,|\,\g_i\le\a_i,\,\forall\,i\in I \}.
\eqno(2.8)$$ Denote $ ({}^\a_\g)=\prod_{i\in I
}({}^{\a_i}_{\g_i}),\
  \forall\,\a,\g\in J ,\,x\in H(A).
$ Then $({}^\a_\g)=0$ if $\g\not\in J (\a)$. We have
$$
\p^{\a}(xy)=\sum_{\g\in{
J}}({}^\a_\g)\e^+({\a-\g,\g})^{-1}\e({{\p^\g},{x}})
\p^{\a-\g}(x)\p^{\g}(y),\ \ \forall\,\a\in J ,\,x,y\in H(A).
\eqno(2.9)$$ One may easily check that
$$
(u\p^{\a})\cdot(v\p^{\b})=\sum_{\l\in{
J}}({}^\a_\l)\e^+({\a-\l,\l})^{-1}
\e({{\p^\l},{v}})\e^+(\l,\b)u\p^{\a-\l}(v)\p^{\b+\l},\,
\forall\,u,v\in H(A),\,\a,\b\in J . \eqno(2.10)$$ Clearly that
formula (2.10) defines a $\G$-graded associative algebra
$(A[D],\cdot)$ such that $A$ is a left $A[D]$-module via (2.5).
For any
$$
x=\sum_{\a\in J }u_\a\p^{\a}\in A[D],\ u_\a\in A, \eqno(2.11)$$
the expression in (2.11) is sometimes not unique. An expression of
$x$ in (2.11) is called {\it principal}, if the integer ${\rm
max}\{|\a|\,\ \Big|\ \, a_{\a}\ne0\}$ is minimal. We denote this
integer by $h(x)$. Set
$$
F_1=A^D=\{u\in A\,|\, D(u)=0\}. \eqno(2.12)$$ From Lemma 2.1 in
[P2], any nonzero $a\in H(F_1)$ is invertible (i.e., $F_1$ is a
graded field) when $A$ is {\it $\G$-graded $D$-simple}. In this
case $(F_1)_-=0$.
\par
{\bf Lemma 2.1}. {\it i) $A\cap A[D]D=0$. ii) The $\e$-center $Z
(A[D])$ of $A[D]$ is $F_1$.}
\par
{\it Proof}. Note that we have assume that $A$ is $\G$-graded
$(\e,\G)$-commutative. Part i) follows from the action of $x\in
A\cap A[D]D$ on $1$.
\par
ii). For all $u\in H(F_1)$ and $\d\in H(D)$, we have
$[\d,u]=\d(u)=0$. So $F_1\subset Z(A[D])$.
Suppose $x\in H(Z (A[D]))$. If $h(x)>0$, write $x=x_0+x_1$ where
$x_0\in H(A)$ and $x_1\in H(A[D]D)$. Then $x_1\ne 0$. Choose $a\in H(A)$
such that $x_1(a)\ne 0$. It follows that
$0=[x,a]=[x_1,a]= x_1(a)+y,
$
where $y\in A[D]D$. From i) we deduce that $x_1(a)=y=0$, a
contradiction. So $h(x)=0$, i.e., $x\in A$. We have
$$
x\p=\e(\bar{x},\bar{\p})\p x=\e(\bar{x},\bar{\p})\p(x)+x\p,\
\forall\,\p\in H(D),
$$
to give $\p(x)=0$ for all $\p\in H(D)$, i.e., $x\in F_1$.
Therefore Lemma 2.1 follows.
\hfill$\Box$
\par
{\bf Theorem 2.2}. {\it Suppose that $A$ is a $\G$-graded
$(\e,\G)$-commutative $F$-algebra with an identity, and $D\subset
Der^\e(A)$ is a color commutative subspace. Then the $\G$-graded
associative algebra $A[D]$ is graded simple if and only if $A$ is
$\G$-graded $D$-simple.}
\par
{\it Proof}. "$\Rar$'': Suppose $A$ is not graded $D$-simple.
Choose a nonzero proper graded $D$-stable ideal $\cal K$. Then
clearly ${\cal K}[D]$ is a nonzero graded color ideal of $A[D]$.
Since $A[D]$ is graded simple, then ${\cal K}[D]=A[D]$, in
particular, $A\subset {\cal K}[D]$. From lemma 2.1 i), we know
that $A\subset\cal K$, a contradiction. Thus $A$ is graded
$D$-simple.
\par
``$\Lar$'': Suppose $L$ is a nonzero graded ideal of $A[D]$. It
suffices to show that $L=A[D]$.
\par
Suppose that $L\cap A=0$. Choose $x\in H(L)\bs\{0\}$ such that
$h(x)$ is minimal. So $h(x)>0$. Write $x=x_0+x_1$ where $x_0\in
H(A)$ and $x_1\in H(A[D]D)$. Then $x_1\ne0$. Choose $a\in H(A)$
such that $x_1(a)\ne0$. It follows that
$$
x'=[x,a]=[x_1,a]=x_1(a)+y\in L,
$$
where $y\in H(A[D]D)$. From (1.4) and the computation (2.10), we
deduce that $h(x')=h(y)<h(x)$. By the minimality of $h(x)$ we
deduce that $x'=0$. Applying Lemma 2.2 i) gives $x_1(a)=y=0$, a
contradiction. Thus $L\cap A\ne0$. It is clear that $L\cap A$ is a
graded $D$-ideal of $A$. Since $A$ is graded $D$-simple, $L\cap
A=A$. In particular $1\in L$. Therefore $L=A[D]$. \hfill$\Box$
\vskip 10pt \centerline{\bf \S3  Simple Lie color algebras of Weyl
type}
\par
In this section we  assume that $A$ is $\G$-graded $D$-simple and
$(\e,\G)$-commutative associative algebra with an identity, the
field $F$ is of characteristic not 2, and $D\subset Der^\e(A)$ is
color commutative. In this section we shall first study the
structure of $A[D]$ as a left $A$-module, then investigate the Lie
structure of $A[D]$. We still use the notation in Section 2.
\par
{\bf Lemma 3.1}. {\it i) If $xd^\a=0$ for some $x\in H(A)\bs\{0\}$
and some $ \a\in J $, then $d^\a=0$. \vs{-3pt}\par ii)
$A\subset[A[D]  ,A[D] ]$.}
\par
{\it Proof}.
i) Let $C=\{x\in A\,|\,x\p^{\a}=0\}$. It is easy to
verify that $C$ is a nonzero $\G$-graded
$D$-ideal of $A$. Since $A$ is $\G$-graded $D$-simple, then $C=A$.
It yields that $\p^{\a}=0$.
\par
ii) From $[x\p,y]=x\p(y)\in[A[D],A[D]]$ for all $x,y\in H(A)$ and
$\p\in H(D)$, we obtain that $AD(A)\subset[A[D] ,A[D]]$. Since
$AD(A)$ is a nonzero $\G$-graded $D$-ideal of $A$ and $A$ is
$\G$-graded $D$-simple, $A=AD(A)$. Thus $A\subset[A[D] ,A[D]]$.
\hfill$\Box$
\par
If ${\rm char\sc\,}F=p>0$, observe that for any $\p\in{\rm
Der}^\e( A)$, one has $\p^p\in{\rm Der}^\e(A)$. For convenience,
no matter whether ${\rm char\sc\,}F=p>0$ or ${\rm char\sc\,}F=0$,
we denote
$$
\DD={\rm Der}^\e(A)\cap F_1[ D], \eqno(3.1)$$ where $F_1[D]$ is
the subalgebra of ${\rm End}(A)$ generated by $F_1$ and $D$. Let
$\{d_i\,|\,i\in {\cal I}\}$ be a homogeneous $F_1$-basis for
$\DD$, where ${\cal I}$ is some index set. It is clear that
$D\subset\DD$, and that $A$ is graded $D$-simple if and only if
$A$ is graded $\DD$-simple. Let
$$
\matrix{ {\cal J}=\!\!\!\!& \{\a=(\a_i\,|\,i\in {\cal I})\ |\
\a_i\in\Z_+,\mbox{ and }\a_i\le p-1\mbox{ if } {\rm
char\sc\,}F=p>0,\, \vs{4pt}\hfill\cr& \mbox{ \ and }\a_i=0\mbox{
or }1\mbox{ if }\bar d_i\in\G_-, \mbox{ and }\a_i=0\mbox{ for all
but a finite number of }i\in {\cal I}\}. \hfill\cr} \eqno(3.2)$$
We will often simply denote $\a=(\a_i\,|\,i\in {\cal I})$ by
$\a=(\a_i)$. We also fix a total ordering $<$ on ${\cal I}$ and
define a total ordering on ${\cal J}$ as in (2.7). We also write
$d^\a=\prod_{i\in I}d^{\a_i}$ according to the ordering $<$ on
${\cal I}$ for $\a\in {\cal J}$. Then
$$
A[D]=\sum_{\a\in {\cal J}}Ad^\a. \eqno(3.3)$$

 {\bf Theorem 3.2}. {\it Suppose that the field $F$ is of
characteristic not 2, $A$ is a $\G$-graded $D$-simple and
$(\e,\G)$-commutative associative $F$-algebra with an identity,
where $D\subset Der^\e(A)$ is a color commutative subspace. Then
$A[D]$ is a free $A$-module with the homogeneous basis
$\{d^\a\,|\,\a\in {\cal J}\}$.}
\par
We break the proof of this theorem into several Lemmas. For ${\cal
J}_0\subset{\cal J},\,B\subset A$, we say that the sum
$\sum_{\a\in {\cal J}_0}Bd^\a$ is {\it direct} if a finite sum
$\sum_{\a\in {\cal J}_0}b_\a d^\a=0$ with $b_\a\in B$ implies
$b_\a d^\a=0$ for all $\a\in {\cal J}_0$.
\par
{\bf Lemma 3.3}. {\it Let ${\cal J}_0\subset {\cal J}$. The sum
$\sum_{\a\in {\cal J}_0}Ad^{\a} $ is direct if and only if the sum
$\sum_{\a\in {\cal J}_0}F_1d^{\a}$ is direct.}
\par
{\it Proof}.
``$\Rar$'': This direction is clear.
\par
``$\Lar$'': Suppose $\sum_{\a\in {\cal J}_0}Ad^{\a}$ is not
direct. There exists ${\cal
J}_1=\{\a^{(0)},\a^{(1)},\cdots,\a^{(r)}\}\subset {\cal J}_0$ with
$|{\cal J}_1|>1$ such that
$Ad^{\a^{(0)}}\cap\sum_{i=1}^rAd^{\a^{(i)}}\ne0$. Choose such a
minimal $r$. Let
$$
ud^{\a^{(0)}} =\sum_{i=1}^ru_id^{\a^{(i)}}\in
(Ad^{\a^{(0)}}\cap\sum_{i=1}^rAd^{\a^{(i)}})\bs\{0\},
\eqno(3.4)
$$
for some $u,\,u_i\in H(A)$, and let $I_1=\{0,1,\cdots,r\}$. It
follows that $d^{\a^{(i)}}\ne0$ for any $i\in I_1$. Let
$$
B_0=\{x\in A\,|\,xd^{\a^{(0)}}\in\sum_{i\in
I_1\bs\{0\}}Ad^{\a^{(i)}}\}.
$$
Then it is easy to see that $B_0$ is a nonzero $\G$-graded
$D$-ideal of $A$. Since $A$ is $\G$-graded $D$-simple, then
$B_0=A$. Then
$$
d^{\a^{(0)}}=\sum_{i\in I_1\bs\{0\}}a_id^{\a^{(i)}}
\mbox{ for some }a_i\in A.
\eqno(3.5)$$
By taking bracket with $\p\in H(\DD)$, we deduce that
$$
0=[\p,d^{\a^{(0)}}]=[\p,\sum_{i\in I_1\bs\{0\}}a_id^{\a^{(i)}}]
=\sum_{i\in I_1\bs\{0\}}\p(a_i)d^{\a^{(i)}}.
$$
From the minimality of $r$, we must have
$$
\p(a_i)d^{\a^{(i)}}=0,\,\forall\,\p\in\DD,\,1\le i\le r.
\eqno(3.6)$$ Suppose $D(a_i)\ne0$ for some $i\in I_1$, say
$D(a_1)\ne0$. Lemma 3.1 yields that $d^{\a^{(1)}}=0$, a
contradiction. Thus $D(a_i)=0$ for all $i\in I_1$, i.e., $a_i\in
F_1$ in (3.5). Therefore the sum $\sum_{\a\in {\cal
J}_1}F_1d^{\a}$ is not direct. This proves the lemma. \hfill$\Box$
\par
{\bf Lemma 3.4}. {\it $A[D]=\oplus_{\a\in {\cal J}}Ad^\a$.}
\par
{\it Proof}. Suppose the sum in (3.3) is not direct. From Lemmas
2.1 i) and 3.1, there exists ${\cal
J}_1=\{\a^{(0)},\a^{(1)},\cdots,\a^{(r)}\}\subset {\cal
J}\bs\{0\}$ such that $F_1d^{\a^{(0)}}\cap\sum_{i=1}^rF_1
d^{\a^{(i)}}\ne0$. Let $I_1=\{0,1,\cdots,r\}$, and let $h({\cal
J}_1)={\rm max}\{|\a^{(i)}|\,\,|\,\,i\in I_1\}$. Then $h({\cal
J}_1)>1$. We choose ${\cal J}_1$ such that $h({\cal J}_1)$ is
minimal, and then $r$ is minimal. Let
$$
-d^{\a^{(0)}}=\sum_{i=1}^ra_id^{\a^{(i)}}\in
(F_1d^{\a^{(0)}}\cap\sum_{i=1}^rF_1d^{\a^{(i)}}) \bs\{0\},
\eqno(3.7)$$ for some $a_i\in H(F_1)$, $1\le i\le r$. Denote
$a_0=1$. It follows that $a_id^{\a^{(i)}}\ne0$ for $i\in I_1$.
Rewrite (3.7) to give $\sum_{i=0}^ra_id^{\a^{(i)}}=0. $ Then for
$x\in H(A),\,\sum_{i=0}^ra_id^{\a^{(i)}}(x)=0$. By (2.10),
\begin{eqnarray*}
0&=&[\sum_{i=0}^ra_id^{\a^{(i)}},x] =\sum_{i=0}^r\sum_{\g\in
J}a_i(^{\a^{(i)}}_{\,\,\g})
\e^+({\a^{(i)}-\g,\g})^{-1}
\e({d^{\g}},{x})d^{\a^{(i)}-\g}(x)d^{\g} \\
&=&\sum_{\g\in J,\,0<|\g|<h(J_1)}\left(
\sum_{i=0}^ra_i(^{\a^{(i)}}_{\,\,\g})\e^+({\a^{(i)}-\g,\g})^{-1}
\e({d^{\g}},{x}) d^{\a^{(i)}-\g}(x)\right)d^{\g}.
\end{eqnarray*}
From the minimality of $h({\cal J}_1)$ and Lemma 2.1 i), it
follows that
$$
\sum_{i=0}^r a_i(^{\a^{(i)}}_{\,\,\g})\e^+({\a^{(i)}-\g,\g})^{-1}
\e({d^{\g}},{x}) d^{\a^{(i)}-\g}(x)d^{\g}=0, \eqno(3.8)
$$
for all $x\in H(A)$ and all $\g\in {\cal J}\setminus\{0\}$ with
$0<|\g|<h({\cal J}_1)$. Noting that for any $\g\in {\cal J}$ with
$\a^{(i)}-\g\in {\cal J}$, we have $d^\g
d^{\a^{(i)}-\g}=\e({d^\g},{d^{\a^{(i)}-\g}}) d^{\a^{(i)}}\ne0$, to
give $d^\g\ne0$. Since $d^{\g}\ne0$, from (3.8) and Lemma 3.1 we
obtain $\sum_{i=0}^r
a_i(^{\a^{(i)}}_{\,\,\g})\e^+({\a^{(i)}-\g,\g})^{-1}
\e({d^{\g}},{x})
 d^{\a^{(i)}- \g}(x)=0,
$
i.e.,
$$
\sum_{i=0}^r a_i(^{\a^{(i)}}_{\,\,\g})\e^+({\a^{(i)}-\g,\g})^{-1}
 d^{\a^{(i)}- \g}(\e({d^{\g}},{x})
x)=0,
$$
for all $x\in H(A)$ and all $\g\in {\cal J}$ with  with
$\a^{(i)}-\g\in {\cal J}$. So $\sum_{i=0}^r a_i
(^{\a^{(i)}}_{\,\,\g})\e^+({\a^{(i)}-\g,\g})^{-1}
d^{\a^{(i)}-\g}=0$ for all $\g\in {\cal J}$ with $0<h(\g)<h({\cal
J}_1)$. We have $a_i(^{\a^{(i)}}_{\,\,\g})\ne0$ since
$\a_j^{(i)}\le p-1$ if ${\rm char\sc\,}F=p>0$, for $\g\in {\cal
J}$ with $\a^{(i)}-\g\in {\cal J}$. Thus we see that
$d^{\a^{(i)}-\g}=0$ for all $\g\in {\cal J}\setminus\{0\}$ with
$\a^{(i)}-\g\in {\cal J}\setminus\{0\}$. As we noticed also
$d^{\a^{(i)}-\g}\ne0$, so this contradicts the minimality of
$h({\cal J}_1)$. Therefore $A[D] =\oplus_{\a\in {\cal J}_0}Ad^\a.$
\hfill$\Box$
\par
{\bf Lemma 3.5}. {\it For all $\a\in {\cal J}$, $d^\a\ne0$.}
\par
{\it Proof}. It is obvious that $d^{\a}\ne0$ for $\a\in {\cal J}$
with $|\a|\le1$. Suppose $d^{\a}=0$ for some $\a\in {\cal J}$ with
$|\a|\ge2$. Choose $\a$ such that $|\a|$ is minimal. Then for any
$x\in H(A)$, $d^{\a}(x)=0$. Then by (1.4), (2.10),
$$
0=[d^{\a},x]=\sum_{\b\in {\cal J}}(^\a_\b) \e^+({\a-\b,\b})^{-1}
\e(\b,\bar{x}) d^{\a-\b}(x)d^{\b},\,\forall\,x\in H(A).
$$
 Applying
Lemma 3.4 and noting that $(^\a_\b)\e^+({\a-\b,\b})^{-1}
\e(\b,\bar{x})\ne0$ (by noting that $\a_i\le p-1$ if ${\rm
char\sc\,}F=p>0$) for any $\b\in {\cal J}$ with $\a-\b\in {\cal
J}$, we get $d^{\a-\b}(x)d^{\b}=0,\,\forall\,x\in H(A), \b\in J$.
Noting that $d^{\b}\ne0$ for any $\b\in {\cal J}$ with $\a-\g\in
{\cal J}\setminus\{0\}$ (since $|\a|$ is minimal), applying Lemma
3.1 gives $d^{\a-\b}(x)=0,\, \forall\,x\in H(A)$, i.e.,
$d^{\a-\b}=0$ for $\b\in {\cal J}$ with $\a-\g\in {\cal
J}\setminus\{0\}$, contradicting the minimality of $|\a|$.
Therefore the assertion follows. \hfill$\Box$
\par
{\it Proof of Theorem 3.2}.
It follows from Lemmas 3.1, 3.3, 3.4 and 3.5.
\hfill$\Box$
\par
Denote by $W=[A[D] ,A[D] ]$ the derived Lie color ideal of $A[D]
$. In the rest of this section We will investigate the structure
of $W$.
\par
{\bf Lemma 3.6}.
{\it Suppose ${\rm dim}_{F_1}\DD>1$ or $\DD\ne\DD_-$. Then for any
homogeneous color-derivation $\p\in\DD$, we have $A\p\subset W$.}
\par
{\it Proof}.
\ul{{\it Case 1}: $\p\in\DD_+$}.
For any $x,a\in H(A)$, by Lemma 3.1 ii), $a,\,x\p^2(a)\in W$, thus
$$
x\e(\p, a)\p(a)\p=\frac{1}{2}([x\p^2,a]-x\p^2(a))\in W,
\eqno(3.9)$$ i.e., $A\p(A)\p\subset W$. But the space $K=\{x\in
A\,|\,x\p\in W\}$ is $\G$-graded $D$-stable since $W$ is a Lie
color ideal of $A[D]$, and $K$ contains the nonzero graded ideal
$A\p(A)$ of $A$, by Lemma 3.1 in [P2], $K=A$, i.e., $A\p\subset
W$.
\par
\ul{{\it Case 2}: $\p\in\DD_-$}. First assume that $\DD_+\ne0$.
Choose $\p'\in H(\DD_+)\bs\{0\}$, by Case 1, $A\p'\subset W$. Thus
for any $x,y\in H(A)$, we have $[x\p',y\p]\in W$ and
$$
[x\p',y\p]=x\p'(y)\p-\e({x\p'},{y\p})y\p(x)\p'\equiv
x\p'(y)\p\mbox{ \ (mod\,}W), \eqno(3.10)$$ i.e., $A\p'(A)\p\subset
W$, thus as in Case 1, $A\p\subset W$. Next assume that $\DD_+=0$,
then by the assumption of the lemma, we can choose $\p'\in\DD_-$
such that $\p,\p'$ are $F_1$-linear independent. For $x,y\in
H(A)$, we have $[x\p',y\p]\in W$ and
$$
\matrix{ [x\p',y\p]\!\!\!\!&=x\p'(y)\p-\e({x\p'},{y\p})y\p(x)\p'
\vs{4pt}\hfill\cr& =x\p'(y)\p-\e({x\p'},{y\p})\e(\p,
y)^{-1}(\p(yx)\p'-\p(y)x\p') \vs{4pt}\hfill\cr& \equiv
x(\p'(y)\p+\e({x\p'},{y\p})\e( y,\p) \e({\p(y)}, x)\p(y)\p')\mbox{
\ (mod\,}W) \vs{4pt}\hfill\cr& =x(\p'(y)\p+\e({\p'},{y\p})\e(
y,\p)\p(y)\p'), \hfill\cr} \eqno(3.11)$$ where the equality
``$\equiv$'' follows from that $\p(xy)\p'=[\p,xy\p']\in W$. Denote
$d=\p'(y)\p+\e({\p'},{y\p})\e( y,\p)\p(y)\p'.$ Note that (3.11)
shows that $Ad\in W$. By using $\p^2=0$ and  applying ${\rm
ad\sc\,}\p$ to (3.11), we obtain $[\p,xd]\in W$ and
$$
[\p,xd]=\p(x)d+\e(\p, x)x[\p,d]\equiv \e(\p, x)x\p(\p'(y))\p\mbox{
\ (mod\,}W).
$$
This shows that $A\p(\p'(A))\p\!\subset\! W$. But
$\p\p'\!\ne\!0$ by Lemma 3.5,
thus as in Case 1, $A\p\!\subset\! W$.
\hfill$\Box$
\par
For any $i\in {\cal I}$, we define
$$
\d^{(i)}\in {\cal J}\mbox{ such that }\d^{(i)}_j=\d_{i,j},\
\forall\,j\in {\cal I}. \eqno(3.12)$$ The following technical
lemma plays a crucial role in describing the structure of $W$.
\par
{\bf Lemma 3.7}. {\it Let $\b\in {\cal J}$ with ${\it
supp}(\b)=\{1,2,...,n\}$ and let $\p'_1,\p'_2,...,\p'_n\in{\rm
Der}^\e( A)$ be $A$-linear independent homogeneous derivations
such that $\bar\p'_i-\bar d_i=\bar\p'_j-\bar d_j$ for
$i,j=1,2,...,n$. Suppose there exist $a_1,a_2,...,a_n\in
F\bs\{0\}$ such that
$$
x\sum_{i=1}^n a_i\e({d^{\b-\d^{(i)}}},b)
\p'_i(b)d^{\b-\d^{(i)}}\in W, \ \forall\,x,b\in H(A),
\eqno(3.13)$$ then $Ad^{\b-\d^{(i)}}\subset W$ for $i=1,2,...,n$.}
\par
{\it Proof}. By shifting the index, it suffices to prove
$Ad^{\b-\d^{(n)}}\subset W$. We shall prove by induction on $n$.
If $n=1$, (3.13) shows that $A\p'_1(A)d^{\b-\d^{(n)}}\subset W$,
thus the result follows from [P2]. Suppose $n\ge2$. Replacing $x$
by $x\p'_1(a)$, and replacing $x,b$ by $y\p'_1(b),a$ in (3.13), we
obtain respectively
$$\matrix{\dis
x\sum_{i=1}^n a_i\e({d^{\b-\d^{(i)}}}, b)
\p'_1(a)\p'_i(b)d^{\b-\d^{(i)}}\in W, \vs{8pt}\hfill\cr\dis
y\sum_{i=1}^n a_i\e({d^{\b-\d^{(i)}}}, a)
\e({\p'_1(b)},{\p'_i(a)}) \p'_i(a)\p'_1(b)d^{\b-\d^{(i)}}
 \in W.
\hfill\cr} \eqno(3.14)$$ Setting $y=x\e({d^{\b-\d^{(1)}}}, a)^{-1}
\e({d^{\b-\d^{(1)}}}, b) \e({\p'_1(b)},{\p'_1(a)})^{-1}$ and
subtracting the two expressions of (3.14), we obtain
$$
x\sum_{i=2}^n a_i\e({d^{\b-\d^{(i)}}}, b)
(\p'_1(a)\p'_i(b)-u_i(a)\p'_i(a)\p'_1(b))d^{\b-\d^{(i)}}\in W,
\eqno(3.15)$$ where
$$
\matrix{ u_i(a)\!\!\!&=\e({d^{\b-\d^{(1)}}}, a)^{-1}
\e({d^{\b-\d^{(1)}}}, b)
\e({\p'_1(b)},{\p'_1(a)})^{-1}\cdot\vs{4pt}\hfill\cr& \hskip
2cm\cdot\e({d^{\b-\d^{(i)}}}, a)^{-1} \e({d^{\b-\d^{(i)}}}, b)
\e({\p'_1(b)},{\p'_i(a)}) \vs{4pt}\hfill\cr&
=\e({d^{\d^{(1)}-\d^{(i)}}}, a) \e(\p'_1,\p'_i\p'_1)\in F.
\hfill\cr} \eqno(3.16)$$ Fix $a\in H(A)$ such that $\p'_1(a)\ne0$,
then $\p''_i=\p'_1(a)\p'_i-u_i(a)\p'_i(a)\p'_1\in{\rm Der}_{F}^\e(
A) ,\,i=2,...,n$ are $A$-linear independent derivations satisfying
the conditions of the Lemma, thus by induction
$Ad^{\b-\d^{(n)}}\subset W$. \hfill$\Box$

\par
{\bf Lemma 3.8} [ZS, Theorem 2.12].
{\it Suppose that $A$ is a $\G$-graded simple associative algebra of
characteristic not 2 and ${\rm dim}_{Z (A)}A>4$. Then
$[A,A]/([A,A]\cap Z )$ is a simple $\e$-Lie color algebra.}
\hfill$\Box$
\par
Our second main result in this section is the following.
\par
{\bf Theorem 3.9}. {\it Suppose that  $F$ is of characteristic not
2, $A$ is a $\G$-graded $D$-simple $(\e,\G)$-commutative
associative $F$-algebra, where $D\subset Der^\e(A)$ is color
commutative and nonzero. Let $W=[A[D] ,A[D] ]$. \vs{-3pt}\par i)
If $|{\cal J}|=\infty$, then $W=A[D] $. \vs{-3pt}\par ii) Suppose
$|{\cal J}|<\infty$. Let $\g\in {\cal J}$ be the maximal element
of ${\cal J}$ (i.e., $|\g|=h({\cal J})$), then $W=\oplus_{\a\in
{\cal J}\bs\{\g\}}Ad^\a\oplus\DD(A)d^\g$. \vs{-3pt}\par iii) The
Lie color algebra $\ol{A[D]}=W/F_1$ is simple except when
$A=F_1[t]/(t^2)$ and $D=F{d\over dt}$ and both $t$ and $d\over dt$
have colors in $\G_-$. }
\par
{\it Proof}. If dim$\DD=1$ and $\DD=\DD_-$, we see that $|{\cal
J}|=2$. Clearly ii) is true in this case. Now suppose that
$\DD\ne\DD_-$ if dim$\DD=1$. Then $|{\cal J}|>2$.

For any $\b\in {\cal J}$ with $|\b|\ge2$. Using inductive
assumption, we may suppose that $Ad^\a\subset W$ for all $\a\in
{\cal J}(\b)$ with $|\a|\le|\b|-2$. Assume that ${\rm
supp}(\b)=\{1,2,...,n\}$. Then for all $x,b\in H(A)$ we have
$[xd^\b,b]\in W$ and
$$
[xd^\b,b]=x[d^\b,a]\equiv x\sum_{i=1}^n\b_i \e({d^{\b-\d^{(1)}}},
b) \e^+(\d^{(i)},\b-\d^{(i)})^{-1} d_i(b)d^{\b-\d^{(i)}} \mbox{
\,(mod\,}W). \eqno(3.17)$$ By Lemma 3.4, $d_1,d_2,...,d_n$ are
$A$-linear independent derivations. By Lemma 3.7,
$Ad^{\b-\d^{(i)}}\subset W$ for $i=1,2,...,n$. This in particular
proves i), and that $\oplus_{\a\in {\cal J}\bs\{\g\}}Ad^\a\subset
W$ if $|{\cal J}|<\infty$. Next we assume that $|{\cal
J}|<\infty$. Then $|{\cal I}|=m<\infty$ and either ${\rm
char\sc\,}F=p>0$ (in this case $\g_i=p-1$ or 1 for $i\in {\cal
I}$) or $\DD=\DD_-$ (in this case $\g_i=1$ for $i\in {\cal I}$).
Then $\DD(A)d^\g=[\DD,Ad^\g]\subset W$. Thus $W'=\oplus_{\a\in
{\cal J}\bs\{\g\}}Ad^\a\oplus\DD(A)d^\g\subset W$. To prove
$W\subset W'$, let $\a,\b\in {\cal J}$ be such that $\a+\b-\g\in
{\cal J}\setminus\{0\}$. It suffices to show the coefficient of
$d^\g$ in $[xd^\a,yd^\b]$ to be in $D(A)$. By (2.7), (2.10), the
coefficient of $d^\g$ in $[ud^\a,vd^\b]$ for $u,v\in H(A)$ is
$$
\matrix{ ({}^{\,\,\a}_{\g-\b}) \e^+(\a+\b-\g,\g-\b)^{-1}
\e({d^{\g-\b}}, v) \e^+({\g-\b,\b})ud^{\a+\b-\g}(v)
\vs{4pt}\hfill\cr \hskip 3cm-\e({ud^\a},{vd^\b})
({}^{\,\,\b}_{\g-\a})\e^+(\a+\b-\g,\g-\a)^{-1} \e({d^{\g-\a}}, u)
\e^+({\g-\a,\a})vd^{\a+\b-\g}(u) \vs{4pt}\hfill\cr
=\e^+(\a+\b-\g,\g-\b)^{-1} \e({d^{\g-\b}},v)\e^+({\g-\b,\b})
\Big(({}^{\,\,\a}_{\g-\b})ud^{\a+\b-\g}(v)\hfill\cr \hskip 5cm
-({}^{\,\,\b}_{\g-\a})\e(\bar
u,\b+\a-\g){\sc\,}d^{\a+\b-\g}(u)v\Big). \hskip 4cm (3.18)}
$$
Noting that char$F=p>0$ or $\g_i=1$, from the condition
$\a,\b,\g-\a-\b\in {\cal J}\setminus\{0\}$ we see that
$$({}^{\,\,\a}_{\g-\b})=({}^{\,\,\a}_{\a+\b-\g})=(-1)^{|\a+\b-\g|}({}={\,\,\b}_{\g-\a}).\eqno(3.19)$$
Applying $d_i$ to $u'v'$ for any $u'v'\in H(A)$, we have
$$
d_i(u')v'=-\e( u', d_1)u'd_i(v')\in D(A)).
$$
Using this and (3.19),
one can easily deduce that the right-hand side
of (3.18) is in $D(A)$.
This proves ii). To prove iii), we use Lemma 3.7. By i) and ii),
$Z =F_1$ and we have
${\rm dim}_{Z (A[D])}A[D]>4$ except when $A=F_1\oplus F_1x$
with $x^2=0$,
and $D=F{\partial\over \partial x}$. This leads to the
iii).
\hfill$\Box$
\par
{\it Example 3.10}. Let $n\ge2$ and let $x_1,x_2,...,x_n$ be $n$
symbols which have colors in $\G_-$. Let $A$ be the free
$\e$-commutative associative algebra generated by
$x_1,x_2,...,x_n$. Clearly $A$ has dimension $2^n$. Let $D$ be the
space spanned by the derivations $\p_i={\p\over\p x_i},
\,i=1,2,...,n$. Then we obtain a simple Lie color algebra
$\ol{A[D]}$ of dimension $2^{2n}-2$. In particular, if we take
$\G=\Z/2\Z$, $\e(i,j)=(-1)^{ij},\,i,j\in\Z/2\Z$, then
$A=\Lambda^n(x_1,x_2,...,x_n)$ is the exterior algebra and we
obtain the simple Hamiltonian Lie superalgebra $\ol{A[D]}=H(2n)$
(see \S3.3 in [K]).\small \vskip 20pt \cl{{\bf References}}

\begin{description}

\item[{[J]}] D.A. Jordan, On the simplicity of Lie algebras of
derivations
   of commutative algebras, {\it J. Alg.}, {\bf228} (2000), 580-585.

\item[{[K]}] V.G. Kac, Lie superalgebras, {\it Adv.  Math.}, {\bf26}
   (1977), 8-96.

\item[{[O]}] J.M. Osborn,  New simple infinite-dimensional Lie algebras
of characteristic $0$, {\it J. Alg.}, 185(1996), no.3, 820-835.

\item[{[P1]}] D.S. Passman, Simple Lie algebras of Witt type,
   {\it J. Alg.}, {\bf 206} (1998), 682-692.

\item[{[P2]}] D.S. Passman, Simple Lie color algebras of Witt type,
   {\it J. Alg.}, {\bf 208} (1998), 698-721.
\par
\item[{[SZ1]}] Y. Su, K. Zhao, Simple algebras of Weyl type,
    {\it Science in China A},  {\bf44}(2001), 419-426.
\par
\item[{[SZ2]}] Y. Su, K. Zhao, Isomorphism classes and
automorphism groups
    of algebras of Weyl type, {\it Science in China A},
    {\bf 31}(2001), 961-972.
\par
\item[{[Z1]}] K. Zhao, Simple algebras of Weyl type, II,
    {\it Proc. Amer. Math. Soci.},  {\bf130}(2002), 1323-1332.
\par
\item[{[ZS]}] K. Zhao, Y. Su,
Simple Lie color algebras from graded associative algebras, to
appear.
\end{description}

\par
Addresses: Department of Mathematics, Shanghai Jiaotong
University, Shanghai 200030,  China. Email:
kfimmi@public1.sta.net.cn.
\par
Academy of
Mathematics and system Sciences, Chinese Academy of Sciences,
Beijing 100080,  China. Email: kzhao@math08.math.ac.cn.
\par
 Department of Mathematics, Nanjing
University, Nanjing 210093,  China. Email: \break
lszhu@pub.sz.jsinfo.net.

\end{document}